\documentclass[12pt, a4paper]{article}
\usepackage{amsmath}
\usepackage{amsfonts}
\usepackage{amssymb}
\usepackage{graphicx}
\usepackage[usenames]{color}
\usepackage{colortbl}
\usepackage[latin1]{inputenc}
\usepackage{hyperref}
\usepackage{esint}
\usepackage{color,soul}

 \newcommand{\R}{\mathbb R}
 \newcommand{\N}{\mathbb N}
 \newcommand{\dis}{\displaystyle}
 \renewcommand{\d}{\mathrm d}
 
 \newcommand{\qed}{\Box}
 
 \newcommand\dlim{\displaystyle \lim}
\newcommand{\beq}{\begin{equation}}
\newcommand{\eeq}{\end{equation}}

\newcommand{\md}{\medskip}
\newcommand{\sm}{\smallskip}
\newcommand{\nd}{\noindent}

\newtheorem{thm}{Theorem}
\newtheorem{lem}[thm]{Lemma}
\newtheorem{prop}[thm]{Proposition}
\newtheorem{cor}[thm]{Corollary}
\newtheorem{rem}{\it Remark}
\newtheorem{exmp}{\it Example}
\newenvironment{pf}{\it Proof. \enskip \rm}{\hfill $\qed$}

\begin{document}

\title{Solvability of a doubly singular boundary value problem arising in  front propagation for reaction-diffusion equations}
\author{ Cristina Marcelli$^{(*)}$}

\date{}

\maketitle
\noindent
\begin{small}$\,^{(*)}$Dipartimento di Ingegneria Industriale e Scienze Matematiche, Universit\`a Politecnica delle Marche, Via Brecce Bianche 12, Ancona I-60131, Italy; 
c.marcelli@univpm.it.
\end{small}

\bigskip

\begin{small}
\noindent {\bf Abstract.}  The paper deals with the solvability of the following doubly singular boundary value problem
\[\begin{cases} \dot z = c g(u)-f(u) -\dfrac{h(u)}{z^\alpha}\\  z(0^+)=0, z(1^-)=0, \ z(u)>0 \text{ in } (0,1)\end{cases}\]
naturally arising in the study of the existence and properties of travelling waves for reaction-diffusion-convection equations governed by the \(p-\)Laplacian operator.

Here \(c,\alpha\) are real parameters, with \(\alpha>0\), and \(f,g,h\) are continuous functions in \([0,1]\), with 
\[ h(0)=h(1), \quad h(u)>0 \text{ in } (0,1).\]

\bigskip

{\bf AMS} Subject  Classifications: 35C07, 35K57, 34B15, 34B16, 34B08.

\medskip

\end{small}

\section{Introduction}

One of the most popular topics in the study of diffusion equations, in their various possible variants (reaction-diffusion, convection-diffusion, aggregation-reaction-diffusion, etc...) is the existence and the properties of travelling fronts. These particular solutions \(v\), satisfying \(v(\tau,x)=u(x-c\tau)\) for some one-variable function \(u\), have a relevant role 
 in understanding asymptotic behavior of the other solutions of the PDE (see, e.g., \cite{mr,DK,DT2}).

When searching for travelling fronts, the original PDE reduces to an autonomous ordinary differential equation and if  the fronts are monotone this can be further reduced to a first-order singular ordinary equation or, equivalently, to a singular integral equation.
For instance, the travelling fronts of the general reaction-diffusion-convection equation with accumulation term 
\[f(v)v_x+g(v) v_\tau = ({\cal D}(v) v_x)_x +\rho(v) \]
are solutions of the following autonomous second-order equation
\[ (D(u)u')' +(cg(u)-f(u))u'+\rho(u)=0\]
where \(c\) is the wave speed and 
\({}^\prime\) means the derivative with respect to the wave variable \(t\). 
Therefore, the existence of monotone fronts is equivalent to the solvability of the following first-order equation
\[ \dot z= cg(u)-f(u)-\frac{D(u)\rho(u)}{z}.\]
Hence, when the reaction term \(\rho\) vanishes at \(0\) and \(1\)  (and it is positive elsewhere), then the study of the existence of monotone travelling fronts, connecting the equilibria \(0\) and \(1\), involves the following doubly singular boundary value problem:
\beq \begin{cases} \label{pr:no lapl} \dot z = cg(u)-f(u)-\frac{h(u)}{z} \\
z(0^+)=z(1^-)=0\\ 
z(u)>0 \ \text{ in } (0,1)
 \end{cases}\eeq
 where we have put \(h(u):=D(u)\rho(u)\) (see \cite{cmp1}).
 
Problem \eqref{pr:no lapl} admits an equivalent formulation as a singular Volterra-type integral equation:
\[ \begin{cases}  z(u)= \dis\int_0^u  \left(cg(u)-f(u)-\frac{h(u)}{z}\right) \d u\\
z(1^-)=0\\ 
z(u)>0 \ \text{ in } (0,1)
 \end{cases}\]
(see \cite{Gil} in the case of constant \(g\)).
 
\medskip
In the last years, diffusion equations governed by the \(p-\)Laplacian operator, that is 
\[f(v)v_x+g(v) v_\tau = \Delta_p(v) +\rho(v) \]
 have been the subject of growing interest (see, e.g. \cite{AV1,GS2,Aud1,AV2,DT1}). In this case, problem \eqref{pr:no lapl} becomes 
 \beq \begin{cases} \label{pr:sing} \dot z = cg(u)-f(u)-\dfrac{h(u)}{z^\alpha} \\
z(0^+)=z(1^-)=0\\ 
z(u)>0 \ \text{ in } (0,1).
 \end{cases}\eeq
for a suitable \(\alpha>0\) depending on \(p\).

Despite a very extensive literature on the subject, a complete study about the solvability of problem \eqref{pr:sing} seems not known, at least in such a general form (se \cite{Ga} for a recent result without the convection and accumulation terms \(f\), \(g\)).

\bigskip
In this paper we provide an existence and non-existence result for problem \eqref{pr:sing}, given in terms of admissible wave speeds \(c\), related to the value of the parameter \(\alpha\).
More in detail, the main result is the following: 
\begin{thm}\label{t:main} 
Let \(f,g,h\in C([0,1])\) be such that 

\beq
\label{ip:h} h(0)=h(1)=0 \quad \text { and } \quad
h(u)>0 \ \text{ for every } u\in (0,1).
\eeq

\nd
Assume that \(g(0)>0\) and \(\dis\int_0^u g(s)  \d s>0\) for every \(u\in [0,1]\). Moreover, assume  that
\beq \text{ there exists the limit } \  h_{\alpha,0}:= \lim_{u\to 0^+}\dfrac{h(u)}{u^\alpha} \in [0,+\infty]. \label{eq:h0}
\eeq 

Then, if \(h_{0,\alpha}=+\infty\) problem \eqref{pr:sing} does not admit solutions for any \(c\in \R\).

\nd
Otherwise, if \(h_{0,\alpha}<+\infty\),
there exists a threshold value \(c^*\) such that problem \eqref{pr:sing} admits solution if and only if \(c\ge c^*\).
Moreover, put
\beq \label{df:const} G_0:= \inf_{u\in (0,1)} \fint_0^u g(s) \ \d s\ , \quad F_0:=\sup_{u\in (0,1)} \fint_0^u f(s) \ \d s  \ , \quad  
H_0:= \sup_{u\in (0,1)} \fint_0^u\frac{h(s)}{s^\alpha} \d s,
\eeq
(where \(\fint\) stands for the mean  value) we have

\beq
\frac{f(0)}{g(0)} + \frac{\alpha+1}{g(0)}\left(\frac{h_{0,\alpha}}{\alpha^\alpha}\right)^\frac{1}{\alpha+1}\le c^*\le \frac{F_0}{G_0}+\frac{\alpha+1}{G_0} \left( \frac{H_0}{\alpha^\alpha}\right)^\frac{1}{\alpha+1}.
\label{eq:stimanuova}
\eeq
Finally, for every \(c\ge c^*\) the solution is unique.

\end{thm}

\bigskip

As a consequence of the previous theorem, we have that problem \eqref{pr:sing} admits a solution for some \(c\in \R\) if and only if the value \(h_{0,\alpha}\) is finite.
 This places an upper bound on the admissible values of \(\alpha\); namely, if \(h(u)\sim Cu^q\) as \(u\to 0^+\), then \(h_{0,\alpha}<+\infty\) if and only if  \( \alpha\le q\).

\bigskip
Finally, note that when \(g\) is constant, say \(g(u)\equiv 1\), then \eqref{eq:stimanuova} reduces to
\[
f(0) + (\alpha+1)\left(\frac{h_{0,\alpha}}{\alpha^\alpha}\right)^\frac{1}{\alpha+1}\le c^*\le \sup_{u\in (0,1)} \fint_0^u f(s) \d s+(\alpha+1)\left( \frac{1}{\alpha^\alpha} \sup_{u\in (0,1)} \fint_0^u\frac{h(s)}{s^\alpha} \d s\right)^\frac{1}{\alpha+1}.\]
When, in addition, \(\alpha=1\) (the standard case) we obtain
\[
f(0) + 2\sqrt{\dot h(0)} \le c^* \le  \sup_{u\in (0,1)} \fint_0^u f(s) \d s+2 \sqrt{   \sup_{u\in (0,1)} \fint_0^u\frac{h(s)}{s} \d s}\]
that is the estimate proved in \cite{mp1}.

\section{Preliminary results}

In this section we establish some preliminary results
concerning the solutions of the equation
\beq \label{e:eq-sing}
\dot z= cg(u)-f(u) -\frac{h(u)}{z^\alpha}
\eeq
for  given \(c\in \R\), \(\alpha>0\), and  continuous functions \(f,g,h\) in \((0,1)\) satisfying \eqref{ip:h}.

\medskip

The first Lemma 
 concerns the maximal existence interval for the solutions of \eqref{e:eq-sing} and the behavior at the extrema of the interval \((0,1)\).

\bigskip

\begin{lem}
\label{l:general} 
Let \(z\in C^1(a,b)\) be a positive solution of equation \eqref{e:eq-sing}, where \((a,b)\subset (0,1)\) is its maximal existence interval. Then, \(a=0\) and both the limits \(z(0^+)\) and \(z(b^-)\) exist and are finite.

Moreover,  if \eqref{eq:h0} holds, then 
  \(h_{\alpha,0}<+\infty\),  and  there exists \(\dot z(0)\). Furthermore, \(\dot z(0)\)  is  a zero of the function  \(\eta_0(t):=
|t|^{\alpha+1} -(cg(0)-f(0))|t|^\alpha+h_{\alpha,0}\).

 Similarly, if \
\(b=1\) and 
\[  \text{ there exists the limit } \ h_{\alpha,1}:= \lim_{u\to 1^-}\dfrac{-h(u)}{(1-u)^\alpha}\in [-\infty,0],\]
then  \(h_{\alpha,1}>-\infty\)  and  there exists \(\dot z(1)\). Moreover,  \(\dot z(1)\)  is a zero of the function \(\eta_1(t):=
|t|^{\alpha+1} +(cg(1)-f(1))|t|^\alpha+h_{\alpha,1}\).

\end{lem}

\begin{pf}
Note that \(\dot z (u) > cg(u)-f(u)\) in \((a,b)\), so  
\(\dot z\) is bounded below, hence, we deduce that  both the limits 
$z(a^+):=\lim_{u\rightarrow a^+}z(u)$ and $z(b^-):=\lim_{u\rightarrow b^-}z(u)$ exist (finite or infinite).
If \(z(a^+)=+\infty\), then equation \eqref{e:eq-sing} implies that
\( \dlim_{u\to a^+} \dot z(u)= cg(a)-f(a)\in \R\),
a contradiction. An analogous reasoning  works for  \(z(b^-)\).  
Therefore, if  \(a>0\), since \((a,b)\)
is the maximal existence interval of \(z\),
we infer \( z(a^+)=0\), implying
 \(\dlim_{u\to a^+}(cg(u)-f(u))z^\alpha(u)-h(u)=-h(a)<0\). Hence, given a real  \(0<\varepsilon< h(a)\), we have 
\( \dot z(u)=cg(u)-f(u)-\dfrac{h(u)}{z^\alpha(u)}< -\dfrac{\varepsilon}{z^\alpha(u)}<0\) in a right neighbourhood of \(a\), which is a contradiction. This implies  \(a=0\). 

\smallskip
 Assume now  \eqref{eq:h0} and put 
\beq
\zeta(u):= \frac{z(u)}{u}, \quad
L:=\limsup_{u\rightarrow 0^+}\zeta(u), \quad \ell:=\liminf_{u\rightarrow 0^+}\zeta(u). \label{eq:position}
\eeq
Suppose, by contradiction,  $\ell<L$.
 So, for every \(k \in (\ell, L)\)  there exist decreasing sequences \((u_n)_n\), \((v_n)_n\), converging to \(0\), such that  
\[\zeta(u_n)=\zeta(v_n)=k, \quad \dot \zeta(u_n)\ge 0, \ \dot\zeta(v_n)\le 0 \quad \text{ for every } n\ge 1.\]

\md\nd
Since \(\dot\zeta(u)= 
\frac{1}{u}(\dot{z}(u)-\zeta(u))
\),
we get
\(
\dot{z}(u_{n})\geq\zeta(u_n)=k
\), 
so
\[ k\le \dot z(u_n)=
cg(u_n)-f(u_n)-\frac{h(u_{n})}{z^\alpha(u_n) }=cg(u_n)-f(u_{n})-\frac{h(u_{n})}{(k u_{n})^\alpha }.
\]
Passing to the limit as $n\rightarrow+\infty$, we have that \(h_{0,\alpha}\) is finite and
\(
k\le cg(0) -f(0)-\frac{h_{\alpha,0}}{k^\alpha},
\)
that is
\[
k^{\alpha +1}-(cg(0)-f(0))k^\alpha+ h_{\alpha,0}\le 0.\]

\noindent	
By means of a similar argument, by using the sequence \((v_n)_n\) instead of \((u_n)_n\), it is possible to show that 
\(
k^{\alpha+1}-(cg(0)-f (0))k^\alpha+ h_{\alpha,0}\ge 0\), hence we conclude 
\[
k^{\alpha+1}-(cg(0)-f(0))k^\alpha+h_{\alpha,0}= 0 \quad \text{ for all } k\in (\ell,L), \]
a contradiction. Then,  \(\ell=L\) and the  limit 
\(\lambda:=\dlim_{u\to 0^+} \frac{z(u)}{u}\ge 0\) exists. 

Finally, since
\beq \dot z(u) =  cg(u)- f(u) -\dfrac{h(u)}{u^\alpha}\cdot \dfrac{u^\alpha}{z^\alpha(u)},\label{e:limit}\eeq
then \(\lambda<+\infty\), otherwise \(\dot z(u)\to cg(0)- f(0)\) as \(u\to 0^+\), a contradiction; so \(\lambda \in \R\).

If  \(\lambda> 0\) then by \eqref{e:limit} we deduce  the existence of the limit \(\dot z(0^+):=\dis\lim_{u\to 0^+} \dot z(u)\) too; moreover,  its value necessarily coincides with  \(\lambda\). So, 
passing again to the limit as \(u\to 0^+\) in \eqref{e:limit}, we infer
\[\lambda^{\alpha+1} = (cg(0)-f(0))\lambda^\alpha - h_{\alpha,0}\]
that is, \(\dot z(0)\) is a zero of the function \(\eta_0(t)\).
Whereas, if \(\lambda=0\), then  \(h_{\alpha,0}=0\) too, otherwise \(\dot z(u)\to +\infty\), a contradiction. Then, even if \(\lambda=0\), it is a zero of the function \(\eta_0(t)\).

\sm
The local analysis near the point \(1\)  can be made by the same way, putting in \eqref{eq:position}
\[ \zeta(u)= \frac{z(u)}{1-u}, \quad L:=\limsup_{u\rightarrow 1^-}\zeta(u), \quad \ell:=\liminf_{u\rightarrow 1^-}\zeta(u).\]
\end{pf}   

\bigskip
In whats follows we will make use of 
the method of lower and upper-solutions. 

Recall that a function \(z\in C^1(I)\), with  \(I\subset (0,1)\), is  a {\em lower-solution} [{\em upper-solution}] for equation \eqref{e:eq-sing} if 
\[
\dot z \le\ [ \ge ]\ cg(u) -f(u) - \dfrac{h(u)}{z^\alpha} \quad  \text{ for all } u\in I.\]

\medskip 
For the reader's convenience, we recall a classic comparison result (see e.g. \cite[Theorems 9.5 - 9.6]{SZ})

\begin{lem}
\label{lem:confronto}
Let \(z\) be a solution of \eqref{e:eq-sing} in  \(I\subset (0,1)\)
and let \(u_0\in I\) be fixed.

Let \(y\) be an upper-solution of \eqref{e:eq-sing} in the same interval, then
\begin{description}
\item{-} \ if \(z(u_0)\le y(u_0)\), we have  \(z(u)\le y(u)\) for all \(u\ge u_0\)
\item{-} \ if \(z(u_0)\ge y(u_0)\), we have  \(z(u)\ge y(u)\) for all \(u\le u_0\).
\end{description}
Whereas, if \(y\) is a lower-solution, then
\begin{description}
\item{-} \ if \(z(u_0)\ge y(u_0)\), then  \(z(u)\ge y(u)\) for all \(u\ge u_0\)
\item{-} \ if \(z(u_0)\le y(u_0)\), then  \(z(u)\le y(u)\) for all \(u\le u_0\).
\end{description}

\end{lem} 

\bigskip

The following result provides 
 the uniform boundedness by below of the solutions of \eqref{e:eq-sing}  on the compact subintervals of \((0,1)\).
 
\bigskip

\begin{lem} Let us fix \(c_0\in \R\) and \(\alpha>0\). \label{lem:equilim} 
For every \(r\in (0,\frac12)\) there exists \(\delta=\delta_r>0\) such that for every \(c\in \R\) with \(|c-c_0|<\delta\) and for positive solution \(z_c\) of equation \eqref{e:eq-sing}, defined in  \((0,1)\), we have 
\[ z_c(u)\ge \delta \quad \text{ for every } u\in [r, 1-r].\]
\end{lem}

\begin{pf}
Let us fix \(r\in (0,\frac12)\) and put  \(m:=\dis\min \{h(u): u\in [r,1-r]\}>0\). 
Let us consider  the map \((c,u,z)\mapsto (cg(u)-f(u))z^\alpha-h(u)\), which is uniformly continuous in the compact set \([c_0-1,c_0+1]\times[0,1]\times [0,m]\). Hence, there exists \(\delta=\delta_r< (mr)^\frac{1}{\alpha+1} \),  such that if \(u_0\in [r,1-r]\), then we have
\beq
(cg(u)-f(u))z^\alpha- h(u) <  -h(u_0) + \frac{m\alpha}{\alpha+1} \le - \frac{m}{\alpha +1}  \text{ when } |u-u_0|,  |z|,  |c-c_0|<\delta. \label{e:psi-delta}
\eeq

For every fixed \(u_0\in [r,1-r]\), define the function
\[\psi(u):=\left(\delta^{\alpha+1} -m (u-u_0)\right)^{\frac{1}{\alpha+1}} \quad \text{  for } u_0\le u\le u_0+ \delta^{\alpha+1}/m,\] where 
the extremum \(u^*:= u_0+ \delta^{\alpha+1}/m < 1\), by the choice of \(\delta\).

Observe that \(\psi(u)\le \delta\) for every \(u\in (u_0,u^*)\), so 
 by \eqref{e:psi-delta} we can deduce that for every \(c\) with \(|c-c_0|<\delta\) we have
 \[\dot \psi(u)> cg(u)-f(u)-\dfrac{h(u)}{\psi^\alpha (u)}, \quad \text{ for every } u\in (u_0, u^*)\] that is \(\psi\) is an upper-solution for equation \eqref{e:eq-sing} in such an interval, with \(\psi(u_0)=\delta\) and \(\psi(u^*)=0\). 
 
Therefore, 
 if \(z_c\) is a solution of problem \eqref{pr:sing} defined in  \((0,1)\), then  \(z_c(u_0)< \delta=\psi(u_0)\). Indeed, since \(\psi\) is an upper-solution for equation  \eqref{e:eq-sing} in the interval \((u_0,u^*)\), then necessarily \(z_c(u)\le \psi (u)\) in the same interval, implying that \(z_c(u^*)=0\), which is a contradiction since \(z_c\) is defined and positive on the whole interval \((0,1)\). Therefore, we conclude that \(z_c(u_0)\ge \delta\). 
 
 The assertion follows from the arbitrariness of \(u_0\).

\end{pf}

\bigskip

\begin{lem} \label{l:converg}  Let \((c_n)_n\) be a decreasing (not necessarily strictly decreasing) sequence converging to \(c_0\) and let 
  \((z_n(u))_n \) be a sequence of positive solutions of equation \eqref{e:eq-sing} in \((0,1)\) for \(c=c_n\), pointwise convergent to a function \(z_0(u)\) in \((0,1)\).  Then we have \(z_0(u)>0\) for every \(u\in (0,1)\) and \(z_0\) is a solution of equation \eqref{e:eq-sing} in \((0,1)\) for \(c=c_0\).
\end{lem}

\begin{pf} 
As a consequence of  Lemma \ref{lem:equilim}, by the arbitrariness of  \(r>0\) we have \(z_0(u)>0\) for every \(u\in (0,1)\). Moreover, again by Lemma \ref{lem:equilim}, for every \(r\in (0,\frac12)\)  we have 
\[  c_ng(u)-f(u)\le c_ng(u)-f(u)-\dfrac{h(u)}{z_n^\alpha(u)} \le  c_ng(u)-f(u)-\dfrac{h(u)}{\delta^\alpha}
\quad \text{ for every } u\in [r,1-r],\]
for every \(n\in \N\). So, by virtue of  the Dominated Convergence Theorem  we get that, for every \(u,u^*\in [r,1-r]\), we have
\begin{align*} z_0(u) -z_0(u^*) = \dlim_{n\to +\infty} (z_n(u)-z_n(u^*))& = \dlim_{n\to \infty} \int_{u^*}^u \left( f(s)-c_ng(s) -\dfrac{h(s)}{z_n^\alpha(s)} \right) {\rm d} s \\ &= \int_{u^*}^u \left( f(s)-c_0g(s) -\dfrac{h(s)}{ z_0^\alpha(s)} \ {\rm d} s \right).
\end{align*}
Hence, \(\dot z_0(u)= f(u) -c_0g(u) -\dfrac{h(u)}{ z_0^\alpha(u)} \) for every \(u\in [r,1-r]\). 
Since \(r>0\) is arbitrary,  we get that \(z_0\) solves equation \eqref{e:eq-sing} in the whole interval \((0,1)\).

\end{pf}

\begin{rem}
\label{rem:conv}
\rm In light of the proof of the previous Lemma, the assertion holds even if the sequence \((z_n)_n\) is defined in a subinterval \((a,b)\subset (0,1)\).
\end{rem}

\bigskip

In the following result a comparison criterium is proved, in order to establish  the existence and non-existence  of the solutions of problem \eqref{pr:sing}. 
\bigskip

\begin{prop}
Suppose \eqref{ip:h} and assume that there exists a lower-solution \(\varphi\) for equation \eqref{e:eq-sing} in the whole interval \((0,1)\), such that  
 \(\varphi(0^+)=0\) and \(\varphi(u)>0\) for every \(u\in (0,1)\).
 
Then, there exists a \(C^1-\)function \(z:(0,1)\to\R\), solution of the singular boundary value problem \eqref{pr:sing}, such that \(0<z(u)<\varphi(u)\) for every \(u\in (0,1)\).
 \label{l:soprasol}
\end{prop}

\begin{pf} 
Notice that \(\dot \varphi(u)\le cg(u)-f(u)\), so   \(\dot \varphi\) is bounded above, implying the existence of the limit \(\varphi(1^-)<+\infty\). 

Let us divide the proof into two cases.

\medskip

 Case 1: \ \(\varphi(1^-)>0\). 

\medskip
\nd
For each \(n\in \N\) let  \(z_n\) be  the solution  of equation \eqref{e:eq-sing}, satisfying the condition \(z(1)=\varphi(1^-)/n\). By applying Lemma \ref{l:general}, we derive that \(z_n\) is defined in the whole interval \((0,1]\). Furthermore, by virtue of Lemma \ref{lem:confronto} we get \(z_n(u)\le \varphi(u)\) for every \(u\in (0,1]\).
Since
  \eqref{e:eq-sing} has a unique solution passing through a given point, we have \(\varphi(u)\ge z_n(u)\ge z_{n+1}(u)>0\) for each \(n\in \N\) and every \(u\in (0,1)\). Put \(\zeta(u):=\dlim_{n\in \N} z_n(u)\), in force of Lemma \ref{l:converg}, we get \(\zeta(u)>0\) in \((0,1)\) and it is a solution of equation \eqref{e:eq-sing}. Moreover, since \(z_n(0^+)=0\) for every \(n\in \N\), we have \(\zeta(0^+)=0\). Finally, by definition, \(z_n(1)\to  0=\zeta(1^-).\)

\medskip

Case 2: \ \(\varphi(1^-)=0\). 

\medskip
\nd
For each \(n\in \N\) let  \(z_n\) be  the  solution  of equation \eqref{e:eq-sing}, satisfying the  condition \(z_n(\frac{n}{n+1})=\varphi(\frac{n}{n+1})\), defined in its maximal existence interval \((a_n, b_n)\).
 By Lemma \ref{l:general} we have \(a_n=0\) for every \(n\in \N\); moreover, by Lemma
 \ref{lem:confronto} we have 
 \(z_n(u)\le \varphi(u)\) for every \(u\in (0,\frac{n}{n+1})\) and \(z_n(u)\ge \varphi(u)\) for every \(u\in(\frac{n}{n+1},b_n)\). Hence, since by Lemma \ref{l:general} there exists  \(z_n(b_n^-)\) and it is  finite,    we  also  infer \(b_n=1\) for every \(n\in \N\).
 
   Furthermore, 
 since \(z_n(\frac{n}{n+1})=\varphi(\frac{n}{n+1})\ge z_{n+1}(\frac{n}{n+1})\), again by the uniqueness of the solution of equation \eqref{e:eq-sing} passing through a point, we also derive that \(z_n(u)=z_{n+1}(u)\) for each \(u\in (0,1)\) or \(z_n(u)>z_{n+1}(u)\) for each \(u\in (0,1)\). Therefore, the sequence  \((z_n)_n\) is decreasing. Set \(\zeta(u)=\dlim_{n\to +\infty} z_n(u)\); by applying Lemma \ref{l:converg} (with \(c_n=c\) for every \(n\)) we have  \(\zeta(u)>0\) for all \(u\in (0,1)\) and \(\zeta\) is a solution of equation \eqref{e:eq-sing}.
 Finally, we  have \(\zeta(u)\le \varphi(u)\) in \( (0,1)\), hence \(\zeta(0^+)=\varphi(0^+)=0\) and \(\zeta(1^-)= \varphi(1^-)=0\). 
   
\end{pf}

\medskip

\begin{cor} \label{c:sopraint}
Let \eqref{ip:h} be satisfied. Assume  that
  there exists a continuous positive function \(\psi:(0,1)\to \R\) such that \(\psi(0^+)=0\) and
\beq
\psi(u) \le \int_0^u \left( cg(s)-f(s) -\dfrac{h(s)}{\psi^\alpha(s)}\right) {\rm d} s \quad \text{ for every } u\in (0,1).
\label{ip:sopraint}
\eeq
Then, the singular boundary value problem \eqref{pr:sing} admits a solution \(z\in C^1(0,1)\).
\end{cor}

\begin{pf}
Put \(\varphi(u):= \dis\int_0^u \left(cg(s)-f(s) -\dfrac{h(s)}{\psi^\alpha(s)}\right) {\rm d} s\). Since \(\psi\) is positive, by \eqref{ip:sopraint} also \(\varphi\) is positive. Moreover, \(\varphi\) is differentiable, with 
\[ \dot \varphi(u) =cg(u)- f(u) - \dfrac{h(u)}{\psi(u)} \le cg(u) - f(u) - \dfrac{h(u)}{\varphi(u)} \quad \text{ for every } u\in (0,1).\]
Finally, since \(\varphi(0^+)= 0\), then \(\varphi\) satisfies all the assumptions of Proposition \ref{l:soprasol}. So, there exists a solution \(z\) of problem \eqref{pr:sing} satisfying \(0<z(u)\le \varphi(u)\) for every \(u\in (0,1)\). 
\end{pf}

\bigskip

We conclude the section by stating a result concerning the behavior of a suitable function, 
which will be used in the proof of Theorem  \ref{t:main}. We omit the proof since it is trivial.

\begin{lem}
\label{lem:M}
Let \(\beta,\gamma \) be fixed with \(\gamma \ge 0\),  and let 
\[
M_{\beta,\gamma}(t):= t^{\alpha+1} - \beta t^\alpha + \gamma \  , \quad t\ge 0. 
\]
Then, \(M_{\beta,\gamma}\) has minimum \(\mu_{a,\beta}\) which has the same sign of the value \( (\alpha+1)\left(\frac{\gamma}{\alpha^\alpha} \right)^\frac{1}{\alpha+1} -\beta\).

\end{lem}

\section{Proof of Theorem \ref{t:main} and some examples}

{\em Proof of Theorem  \ref{t:main}}.

First of all, let us prove that if problem \eqref{pr:sing} admits a solution for some \(c\), then it is unique. To this aim, 
 assume by contradiction that for a fixed \(c\), problem \eqref{pr:sing} admits two different solutions \(z_1,z_2\). Since  the differential equation in \eqref{pr:sing} admits a unique solution passing through  a given point, we get \(z_1(u)\ne z_2(u)\) for every \(u\in (0,1)\). So, we have \(z_1(u)<z_2(u)\) for every \(u\in (0,1)\) (or vice versa) and then 
\begin{multline*} 0= z_1(1^-)-z_1(0^+) = \int_0^1\left( f(u)-cg(u) - \dfrac{h(u)}{z_1^\alpha(u)} \right) {\rm d} u \\ < \int_0^1 \left( f(u)-cg(u) - \dfrac{h(u)}{z_2^\alpha(u)} \right) {\rm d} u = z_2(1^-)-z_2(0^+)=0, \end{multline*}
a contradiction.
\smallskip

Let us now assume \eqref{eq:h0}. Notice that  if \(h_{0,\alpha}=+\infty\) then by Lemma \ref{l:general} problem \eqref{pr:sing} does not admit solutions for any \(c\in \R\). So, from now on we assume \(h_{0,\alpha}<+\infty\).

Let us fix a value 
\[ c> \frac{F_0}{G_0}+\frac{\alpha+1}{G_0} \left( \frac{H_0}{\alpha^\alpha}\right)^\frac{1}{\alpha+1}\]
(see \eqref{df:const}).
Then, if we consider the function \(M_{\beta,\gamma}\), defined in Lemma \ref{lem:M}, for \(\beta:=c G_0-F_0\) and \(\gamma:=H_0\), we have that the minimum of the function \(M_{\beta,\gamma}\) is negative. Let \(L>0\) be such that 
\(M_{\beta,\gamma}(L)<0\).
Hence
\( L^{\alpha+1} < (cG_0-F_0)L^\alpha- H_0\), implying
\[ L < cG_0-F_0- \frac{H_0}{L^\alpha} \le c \fint_0^u g(s)  \d s -\fint_0^u f(s)\d s -\fint_0^u \frac{h(s)}{(Ls)^\alpha} \d s \quad \text{ for every } u\in (0,1].\]
Therefore, put \(\psi(u):=Lu\) we have 
\[  \psi(u) <   \int_0^u \left( cg(s) - f(s)\d s + \frac{h(s)}{\psi^\alpha(s)} \d s \right) \quad \text{ for every } u\in (0,1]\]
and from Corollary \ref{c:sopraint} we deduce that problem \eqref{pr:sing} admits a solution.

\md
Let us consider now a value \(c\) satisfying
 \[c< \frac{f(0)}{g(0)} + \frac{\alpha+1}{g(0)}\left(\frac{h_{0,\alpha}}{\alpha^\alpha}\right)^\frac{1}{\alpha+1}.\] 
Put \(\beta:= cg(0)-f(0)\) and \(\gamma:=h_{0,\alpha}\), in this case   by Lemma \ref{lem:M} we have that 
\(M_{\beta,\gamma}(t)>0\) for every \(t\ge 0\). On the other hand, if a solution \(z\) of problem \eqref{pr:sing} exists, then by Lemma \ref{l:general} \(\dot z(0)\) should be a zero of the function \(M_{\beta,\gamma}\), a contradiction. So, in this case problem \eqref{pr:sing} does not admit solutions.

\md
Finally, put \(\Gamma:=\{c:  \text{ problem } \eqref{pr:sing}  \text{ admits solution}\}\) and let \( c^*:=\inf \Gamma\). In view of what we have just observed, the value \(c^*\) satisfies estimate \eqref{eq:stimanuova}. 

\medskip
In order to show that problem \eqref{pr:sing} admits a solution for \(c=c^*\), 
 consider a decreasing sequence \((c_n)_n\) in \(\Gamma\), converging to \(c^*\)
and let \(z_n\) be the solution of problem \eqref{pr:sing} for \(c=c_n\).
Put 
\beq \label{df:M-m}
M:=\max_{(u,c)\in [0,1]\times [c^*,c_1]} (cg(u)-f(u)), \ \quad  m:=\min_{(u,c)\in [0,1]\times [c^*,c_1]} (cg(u)-f(u)).\
\eeq

\medskip
Let  \((r_n)_n\) be a decreasing sequence  converging to 0, with \(r_1<\frac12\), and put \(I_k:=[r_k, 1-r_k]\). Of course, \(I_k\subset I_{k+1}\) and \((0,1)=\dis \cup_{k\ge 1} I_k \).
 
By Lemma \ref{lem:equilim}, there exists a value \(\delta_{r_1}>0\) such that for every \(n\ge 1\) we have
\[ c_n g(u)-f(u) \ge c_n g(u)-f(u)-\frac{h(u)}{z_n^\alpha(u)} \ge 
c_n g(u)-f(u)-\frac{h(u)}{\delta_{r_1}^\alpha} \ , \quad \text{ for every } u\in I_1.\]
Therefore, by \eqref{df:M-m} we have
\[ m-\frac{h(u)}{\delta_{r_1}^\alpha}\le \dot z_n(u) \le 
M, \quad \text{ for every } u\in I_1\]
from which we deduce the equicontinuity of the sequence of functions \((z_n)_n\)  in \( I_1\).
Moreover, since
\(0< z_n(u)\le M\) in \(I_1\),
then the sequence  \((z_n)_n\) is  equibounded in \( I_1\). Then, we can apply  the Ascoli-Arzel\`a Theorem to derive that there exists  a subsequence \((z_n^{(1)})_n\) which uniformly converges to a certain function \(z_0^{(1)}\) in  \( I_1\).

Let us consider now the interval \( I_2 \supset I_1\). By means of the same reasoning, we infer the existence of a   further subsequence  \((z_n^{(2)})_n\)  which 
uniformly converges to a function \(z_0^{(2)}\) in   \(I_2\). Moreover, we have 
\(z_0^{(2)}(u)=z_0^{(1)}(u)\) for every \(u\in I_1\).

Proceeding by this diagonal argument, one shows that for all \(k\in \N\)  the sequence \((z_n^{(k)})_n\) admits a subsequence \((z_n^{(k+1)})_n\) which uniformly  converges to a function \(z_0^{(k+1)}\) in  \( I_{k+1}\), and  \(z_0^{(k+1)}=z_0^{(k)} (u)\) for all \( u\in I_k\). 

Finally, let us define \(\zeta_0:(0,1)\to \R\), as \(\zeta_0(u)=z_0^{(k)}(u)\) if  \(u\in I_k\). By what we have just observed,  \(\zeta_0\) is a well-defined function. Moreover, in each interval \(I_k\) the function \(\zeta_0\) is uniform limit of the sequence of solutions \((z_n^{(k)})_n\) of problem \eqref{pr:sing} for \(c=c_n\). So, by applying Lemma \ref{l:converg} we deduce 
that \(\zeta_0\) is a solution of problem \eqref{pr:sing} for \(c=c^*\) in the interval \(I_k\). By the arbitrariness of \(k\) we conclude that that \(\zeta_0\) is a positive solution on the whole interval \((0,1)\).

\medskip
Observe now that by the monotonicity of the integral function of \(g\) we get 
\[z_n(u)\le \int_0^u(c_ng (s)-f(s)) {\rm d} s \le \int_\alpha^u (c_1 g(s)-f(s)) {\rm d}s\]  for all \(n\in \N\) and \(u\in (0,1)\). So, also 
\(\zeta_0(u) \le \int_0^u(c_1 g(s)-f(s)) {\rm d} s\), for every \(u\in (0,1)\), implying
that \(\zeta_0(0^+)=0\).

\md
Furthermore, by virtue of Lemma \ref{l:general}, we deduce that there exists the limit 
\(\zeta_0(1^-)\in [0,+\infty)\). So,   all the assumptions of  Proposition \ref{l:soprasol} are satisfied by the function \(\zeta_0\) and we can conclude that  problem \eqref{pr:sing} admis a solution for \(c=c^*\).

\medskip
Finally, let us now show that problem \eqref{pr:sing} admits solution for every \(c>c^*\).  To this aim,
 let \(z_*\) be the solution of problem \eqref{pr:sing} for \(c=c^*\). Again by  by the monotonicity of the integral function of \(g\), for every \(c>c^*\),  we have
\[ z_*(u) =  \int_\alpha^u  \left( c^*g(s) -f(s) - \dfrac{h(s)}{z_*^\alpha(s)} \right)  {\rm d} s \le    \int_\alpha^u  \left(  cg(s) - f(s)- \dfrac{h(s)}{z_*^\alpha(u)}\right)  {\rm d} s \]
So,  \(z^*\) satisfies condition \eqref{ip:sopraint} of Corollary \ref{c:sopraint} and we can infer that   problem \eqref{pr:sing} admits solution also for every \(c>c^*\).

\smallskip

\hfill\(\Box\)

\bigskip
\begin{exmp} \rm
Let us consider equation \eqref{e:eq-sing} with 
\[ f(u)\equiv 0  \ , \quad g(u)=u+1 \ ,  \quad h(u)=u^2(1-u).\]
In this case we have \(h_{0,\alpha}<+\infty\) if and only if \(\alpha\le 2\), so equation \eqref{e:eq-sing} admits solution if and only if \(\alpha \le 2\). Moreover, since  \(F_0=0\) and \(g(0)=G_0=1\), when \(\alpha=2\) we have \(h_{0,\alpha}=H_0=1\), so  estimate \eqref{eq:stimanuova} becomes \(3\sqrt[3]{\frac14} \le c^* \le 3 \sqrt[3]{\frac14}\), that is \(c^*=\frac{3}{\sqrt[3]{4}}\).

\end{exmp}

\begin{exmp} \rm
Let us consider equation \eqref{e:eq-sing} with 
\[ f(u)=u  \ , \quad g(u)=1-u \ ,  \quad h(u)=u(1-u).\]
In this case we have \(h_{0,\alpha}<+\infty\) if and only if \(\alpha\le 1\), so equation \eqref{e:eq-sing} admits solution if and only if \(\alpha \le 1\). Moreover, we have \(F_0=G_0=\frac12\) and \(g(0)=1\). Furthermore,  when \(\alpha=1\) we have \(h_{0,\alpha}=H_0=1\), so   estimate \eqref{eq:stimanuova} for \(\alpha=1\) becomes \(2 \le c^* \le 5 \).

\end{exmp}

\section{Conclusions}

In this paper we have established ad existence result for problem \eqref{pr:sing}, stating that there exist a solution if and only if the value \(h_{0,\alpha}\) is finite (see \eqref{eq:h0}). In this case, there are infinitely many admissible wave speeds, whose minimum value satisfies estimate \eqref{eq:stimanuova}.

The present  result  can be used in the study of travelling waves for reaction-diffusion-advection equations governed by the p-Laplacian operator.

\subsection*{Funding Information}
This article was supported by PRIN 2022 -- Progetti di Ricerca di rilevante Interesse Nazionale, ``Nonlinear differential problems with applications to real phenomena'' (2022ZXZTN2). 

\subsection*{Acknowledgements}
The author is member of the Gruppo Nazionale per l'Analisi Matematica, la Probabilit\`a e le loro Applicazioni (GNAMPA) of the Istituto Nazionale di Alta Matematica (INdAM).

\end{document}